\documentclass[12pt,a4paper,oneside]{amsart}

\ifdefined\SMART
\usepackage[paperwidth=12cm,paperheight=24cm,left=5mm,right=5mm,top=10mm,bottom=5mm]{geometry}
\else
\calclayout
\fi

\usepackage{amssymb,latexsym,amsfonts,textcomp,fancyhdr,calc,graphicx}
\usepackage{pdfpages,amsfonts,amsthm,amssymb,latexsym,graphicx,leftidx,fancyhdr}
\usepackage{amsmath,amstext,amsthm,amssymb,amsxtra}
\usepackage[allcolors=blue,allbordercolors=blue,pdfborderstyle={/S/U/W 1}]{hyperref}
\usepackage[ocgcolorlinks]{ocgx2}
\usepackage{dsfont}
\usepackage[framemethod=tikz]{mdframed}
\usepackage{asymptote}
\usepackage{enumerate}
\usepackage{changepage}

\usepackage[normalem]{ulem}
\usepackage{cite}

\usepackage{txfonts,pxfonts,tikz} 
\usepackage{tgschola}
\usepackage[T1]{fontenc}

\usepackage{mathtools}

\usepackage{bbm}

\theoremstyle{plain}
\newtheorem{theorem}{Theorem}
\newtheorem{lemma}{Lemma}

\theoremstyle{definition}

\newtheorem{case[theorem]}{Case}

\newcommand{\beql}[1]{\begin{equation}\label{#1}}
\newcommand{\eeq}{\end{equation}}

\newcommand{\Abs}[1]{{\left|{#1}\right|}}

\newcommand{\Set}[1]{{\left\{{#1}\right\}}}

\newcommand{\one}{{\bf 1}}

\newcommand{\RR}{{\mathbb R}}

\newcommand{\ZZ}{{\mathbb Z}}
\newcommand{\NN}{{\mathbb N}}
\newcommand{\TT}{{\mathbb T}}

\newcommand{\rank}{{\rm rank\,}}

\newcommand{\vol}{{\rm vol\,}}

\newcounter{rem}
\newcounter{step}
\newcounter{mysec}
\newcounter{mysubsec}[mysec]
\newcommand{\eqn}[1]{\refstepcounter{equation}\label{#1}\eqno(\theequation)}

\begin{document}

\sloppy

\title[Common tiles from fractional tiles]{Common tiles from fractional tiles for two lattices}

\author{Jan Greb\'ik}
\address{\href{https://www.mff.cuni.cz/en/iuuk}{Computer Science Institute, Charles University}, Malostransk\'e n\'am. 25, 118 00 Praha 1, Czechia}
\email{grebikj@gmail.com}
\thanks{JG was supported by an ERC CZ grant LL2608 funded by the Czech Ministry of Education, Youth and Sports from the ERC CZ Program.}

\author{Mihail N. Kolountzakis}
\address{\href{http://math.uoc.gr/en/index.html}{Department of Mathematics and Applied Mathematics}, University of Crete, Voutes Campus, 70013 Heraklion, Greece,\newline and\newline \href{https://ics.forth.gr/}{Institute of Computer Science, Foundation of Research and Technology Hellas}, N. Plastira 100, Vassilika Vouton, 70013 Heraklion, Greece}
\email{kolount@uoc.gr}


\begin{abstract}
Let $L_1, L_2\subseteq\RR^d$ be lattices of the same volume and let $D_1, D_2$ be two bounded fundamental domains for $L_1, L_2$. We show that there is a measurable common tile for $L_1$ and $L_2$ contained in $D_1+D_2$. The main tool is a result of Bowen, Kun and Sabok (2021) which shows how to obtain a measurable perfect matching given a fractional perfect matching in hyperfinite graphings.
\end{abstract}

\date{November 2026}
\maketitle
\tableofcontents

\section{Introduction}\label{s:intro}

A measurable set $E\subseteq\RR^d$ is a \textit{tile} for a lattice $L\subseteq\RR^d$ if
$$
\sum_{\ell\in L}\one_E(x+\ell)=1
\quad\text{for almost every }x\in\RR^d.
$$
Thus the translates of $E$ by $L$ partition $\RR^d$ up to a set of measure zero. Every tile of $L$ has measure equal to the volume of $L$. If we demand the equality above to hold for all $x \in \RR^d$ then we call the set $E$ a \textit{fundamental domain} of $L$.

A classical problem in this area is the Steinhaus tiling problem \cite{moser1981research,sierpinski1958probleme}. In the plane it asks for a set $E\subseteq\RR^2$ which is a fundamental domain for every rotated copy $R_\theta\ZZ^2$ of the integer lattice. There are two versions of the problem. 
\begin{itemize}
\item In the \textit{set-theoretic} version no regularity is required from $E$ and one asks for exactly one point to belong to $E$ from each coset of every $R_\theta\ZZ^2$.
\item In the \textit{measurable} version the set $E$ is required to be Lebesgue measurable and the tiling identities are required only almost everywhere.
\end{itemize}
The set-theoretic problem in the plane was solved by Jackson and Mauldin \cite{jackson2002lattice,jackson2002sets,jackson2003survey}. The measurable problem in the plane remains open. It is known that a measurable Steinhaus set cannot be bounded \cite{beck1989lattice,kolountzakis1996problem}, and stronger restrictions at infinity are known \cite{kolountzakis1999steinhaus}. In dimensions larger than two there are no measurable Steinhaus sets \cite{kolountzakis1999steinhaus,kolountzakis2002steinhaus}. We refer to \cite{kolountzakis2017measurable,kiss2024solutions} and the references there for related forms of the problem.

A finite version asks for a measurable common tile for a finite collection of  lattices
$$
L_1,\ldots,L_n\subseteq\RR^d
$$
of the same volume. In \cite{kolountzakis1997multi} a measurable common fundamental domain was obtained under a direct-sum condition on the dual lattices $L_i^*$. For two lattices no additional condition is needed: Han and Wang \cite{han2001lattice} proved that any two lattices of the same volume have a measurable common tile. Related simultaneous tiling questions for functions are considered in \cite{kolountzakis2022functions}.

The common fundamental domains or tiles in these results are in general unbounded. Recently Grepstad and Kolountzakis \cite{grepstad2026bounded} proved that any two lattices of the same volume have a \textit{bounded} measurable common tile. The present paper imposes a prescribed bound. As was the case in \cite{grepstad2026bounded,etkind2026bounded} the desired measurability is obtained using results for the field of Borel combinatorics.

\begin{theorem}\label{th:main}
Let $L_1,L_2\subseteq\RR^d$ be lattices of the same volume and let $D_i$ be bounded measurable fundamental domains for $L_i$. Then there exists a Lebesgue measurable set $E\subseteq D_1+D_2$ which tiles with both $L_1$ and $L_2$.
\end{theorem}

\subsection{Tiling as perfect matching}\label{ss:matching}

It is not hard to give a proof of Theorem \ref{th:main} using Hall's theorem, but  without measurability for the set $E$. Let us show this proof first.

Let $f=\one_{D_1}*\one_{D_2}$ and assume that the two lattices have volume 1. Clearly $f$ is measurable, nonnegative and tiles at level 1 with both $L_1$ and $L_2$.
\begin{equation}\label{fractional-tiling}
\sum_{\ell_1 \in L_1} f(x+\ell_1) = \sum_{\ell_2 \in L_2} f(x+\ell_2) = 1,
\end{equation}
for almost all $x \in \RR^d$. We also have that
$$
\Set{f > 0} \subseteq D_1 + D_2
$$
is a bounded set.
Write $X_i = \RR^d/L_i$, $i=1,2$. Construct a bipartite graph with sides $X_1, X_2$ where $x_1 \in X_1$ is connected to $x_2 \in X_2$ if there is a point $x \in \Set{f>0}$ such that $x = x_i \bmod L_i$, $i=1,2$. For this connection let the weight of the edge be $f(x)$. Thus the edges of this bipartite graph (there can be multiple edges between two points) are in one-to-one correspondence with the points in $\Set{f>0}$.

The condition \eqref{fractional-tiling} means exactly that the weights we have put on the edges of the graph form a {\em fractional perfect matching}, i.e., all the weights of the edges emanating out of a vertex add up to 1. The conclusion of Theorem \ref{th:main} then means that for this weighted bipartite graph there exists a {\em perfect matching}, i.e.\ a selection of the existing edges which touch every vertex exactly once.

The existence of a fractional perfect matching in turn implies that Hall's condition is valid for any finite subset of $X_1$ or $X_2$: the total number of neighbors of any set of $k$ vertices on the left or $k$ vertices on the right side of the bipartite graph is at least $k$.

Since $\Set{f > 0}$ is bounded the degrees of all vertices in the graph are finite: the edges touching $x_1 \in X_1$ are in one-to-one correspondence with the points $x \in \Set{f>0}$ that are in a coset of $L_1$. (In fact the degrees are uniformly bounded.) For such a bipartite graph (i.e.\ with finite degrees) Hall's theorem holds and therefore there exists a perfect matching. Each of the edges kept in this matching corresponds to a unique point in $\Set{f>0}$. The collection of these points forms the set $E \subseteq \Set{f>0}$ that tiles with both $L_1$ and $L_2$. There is of course no claim of measurability here. (This not-necessarily measurable set $E$ can also be obtained as a special case of a recent result of Speegle \cite{speegle2026fractional}.)

To correct this lack of measurability we use a theorem of Bowen, Kun and Sabok \cite{bowen2021perfect} which guarantees the existence of a perfect matching ``contained'' in a fractional perfect matching, under certain assumptions.

\section{The support graphing}\label{s:frac}

From now on we assume, as we may, that $\vol L_1 = \vol L_2 = 1$ and let
$$
f=\one_{D_1}*\one_{D_2}.
\eqn{eq:convolution}
$$
The function $f$ takes values in $[0, 1]$, is continuous and $\Set{f > 0} \subseteq D_1+D_2$. Then 
$$
\sum_{\ell\in L_i}f(x+\ell)=1,
\qquad i=1,2.
$$
for almost every $x$.

Let
$$
H=L_1\cap L_2,
\qquad r=\rank H,
\qquad k=d-r,
$$
and
$$
M=\RR^d/H,
\qquad
\Lambda_i=L_i/H.
$$
Write $\pi:\RR^d\to M$ for the canonical map.
$M$ is isomorphic to $\RR^k\times\TT^r$. The groups $\Lambda_i$ are discrete and cocompact in $M$, and
$$
\Lambda_1\cap\Lambda_2=\Set{0}.
\eqn{eq:lambda-intersection}
$$
Moreover the group $\Lambda_i$ is isomorphic to $\ZZ^k\times F_i$,
where $F_i$ is a finite abelian group.

Define for $x \in \RR^d$ the $H$-periodic function
$$
F(x)=\sum_{h\in H}f(x+h).
$$
The function $F$ is measurable, $0\leq F\leq1$, and, after decomposing $L_i$ into $H$-cosets,
$$
\sum_{\lambda\in\Lambda_i}F(y+\lambda)=1,
\qquad i=1,2,
\eqn{eq:Ftiling}
$$
for almost every $y\in M$. Put
$$
K=\Set{y\in M:F(y)>0}.
$$
Then $K \subseteq M$ is measurable and relatively compact.

Set $X_i = M/\Lambda_i \cong \RR^d/L_i$, with normalized Haar measure. We define a bipartite graph $\mathcal{B}$ with sides $X_1$ and $X_2$ by putting, for every $p\in K$, an edge between
$$
q_1(p)=p+\Lambda_1\in X_1
\quad\text{and}\quad
q_2(p)=p+\Lambda_2\in X_2
$$
with weight $F(p)$ (by \eqref{eq:lambda-intersection} an edge between two elements of $M$ is unique, if it exists). Equation \eqref{eq:Ftiling} says exactly that these weighted edges form a measurable perfect fractional matching. It is everywhere positive because the edge set is defined by $F(p)>0$.

Recall that a \emph{graphing} is a Borel graph $G=(X,E)$
of bounded degree on a standard probability space $(X,\mu)$
which satisfies the mass-transport identity
$$
\int_A \deg_B(x)\,d\mu(x)
=
\int_B \deg_A(x)\,d\mu(x)
$$
for all measurable sets $A,B\subseteq X$, where
$\deg_B(x)$ is the number of neighbors of $x$ in $B$.
A graphing is \emph{almost everywhere hyperfinite} if,
for every $\varepsilon>0$, there is a measurable set
$Z\subseteq X$ with $\mu(Z)<\varepsilon$ such that,
after discarding a null set, every connected component
of the induced graph on $X\setminus Z$ is finite.

\begin{lemma}\label{lm:graphing}
The bipartite graph $\mathcal{B}$ is a bounded-degree probability-measure-preserving bipartite graphing and is almost everywhere hyperfinite.
\end{lemma}

\begin{proof}
The edges incident to a vertex of $X_i$ are labelled by the
points of one $\Lambda_i$-coset which lie in $K$. Since $K$
is relatively compact and $\Lambda_i$ is discrete, the
degrees are uniformly bounded.

Equip each $X_i$ with its normalized Haar measure $\mu_i$,
and the disjoint union $X_1\sqcup X_2$ with the probability
measure $(\mu_1+\mu_2)/2$. Let $m$ denote Haar measure on
$M$, normalized so that a fundamental domain of either
$\Lambda_i$ has measure $1$. Weil's integration formula \cite[Theorem 2.51]{folland2016course} gives, for every nonnegative measurable function $g$ on $M$,
$$
\int_M g(p)\,dm(p) = \int_{X_i}\sum_{p\in q_i^{-1}(x)}g(p)\,d\mu_i(x).
$$
Consequently, for every nonnegative measurable function
$\Phi$ on the edges of the bipartite graph,
$$
\int_{X_1}\sum_{y\sim x}\Phi(x,y)\,d\mu_1(x)
=
\int_{X_2}\sum_{x\sim y}\Phi(x,y)\,d\mu_2(y).
$$
This is the mass-transport principle \cite[Section 18.4.1]{lovasz2012large}, and shows that the bipartite graph is probability-measure-preserving.

On $X_1$ join two vertices when they have a common neighbor in $X_2$. Call this auxiliary graph $G$. If the corresponding edge labels are $p,p'\in K$, then
$$
p'-p\in\Lambda_2\cap(K-K).
$$
The set on the right is finite.

Recall that if a finitely generated group $\Gamma$ acts by
measure-preserving Borel transformations on a probability
space $(X,\mu)$, and $S$ is a finite symmetric generating
set of $\Gamma$, the associated \emph{Schreier graphing}
has vertex set $X$, with two distinct vertices $x,y$
adjacent whenever $y=s\cdot x$ for some $s\in S$.
This is a bounded-degree graphing.

In our case, the finite set $\Lambda_2\cap(K-K)$ is
symmetric and can be included in a finite symmetric
generating set $S$ of $\Lambda_2$.
Thus $G$ is a Borel subgraph of the Schreier graphing
of the translation action $\Lambda_2\curvearrowright X_1$
associated with $S$.
Since $\Lambda_2$ is finitely generated abelian, it is
amenable, and this Schreier graphing is almost everywhere
hyperfinite \cite[Section 2]{bowen2021perfect}.

To deduce the hyperfiniteness of the original bipartite
graph $\mathcal B$ from that of the auxiliary graph $G$,
let $D$ be a uniform bound for the degrees in $\mathcal B$.
Given $\varepsilon>0$, choose a measurable set
$A\subseteq X_1$ with $\mu_1(A)<\varepsilon/(1+D)$
such that the restriction of $G$ to $X_1\setminus A$
has only finite connected components.

Let $N_{\mathcal B}(A)\subseteq X_2$ be the set of
neighbors of $A$ in $\mathcal B$. Remove from
$\mathcal B$ the vertices in $A\cup N_{\mathcal B}(A)$.
By the mass-transport identity,
\[
\mu_2(N_{\mathcal B}(A))
\leq D\mu_1(A).
\]
Thus the sum of the measures of the removed sets
is less than $\varepsilon$.

We claim that the remaining bipartite graph has only
finite connected components. Indeed, any two neighbors
in $X_1$ of a vertex $y\in X_2$ are adjacent in $G$.
Consequently, all the remaining neighbors of $y$
belong to a single connected component of
$G|_{X_1\setminus A}$. 
Consider a connected component of the remaining bipartite
graph $\mathcal B$, and let $C$ be the set of its vertices
in $X_1$. Any two vertices of $C$ can be joined by a path
in the remaining bipartite graph. Taking only the
$X_1$-vertices of this path gives a walk in
$G|_{X_1\setminus A}$. Hence $C$ is contained in a single
connected component of $G|_{X_1\setminus A}$ and is
therefore finite. Every vertex of the bipartite component
in $X_2$ is a neighbor in $\mathcal B$ of some vertex
of $C$. Since $\mathcal B$ has finite degrees,
there are only finitely many such vertices. Thus every
component of the remaining bipartite graph is finite.
Hence $\mathcal B$ is almost everywhere hyperfinite.

\end{proof}

We shall also use the following observation at the end of the proof.

\begin{lemma}\label{lm:lift}
If the bipartite graph $\mathcal{B}$ has a measurable perfect matching, then there is a measurable set $E\subseteq\Set{f>0}$ which tiles with both $L_1$ and $L_2$.
\end{lemma}

\begin{proof}
Let $\mathcal M$ be a measurable perfect matching. Since there is at most one edge between any pair of vertices, the matching $\mathcal M$ corresponds to a measurable set
$$
E_{\mathcal M}\subseteq K
$$
which meets almost every $\Lambda_i$-coset in exactly one point, for $i=1,2$.

Choose a Borel section $s:M\to\RR^d$ of $\pi:\RR^d\to M$ and enumerate $H=\Set{h_1,h_2,\ldots}$. For $p\in E_{\mathcal M}$ we have
$$
F(p)=\sum_{h\in H}f(s(p)+h)>0.
$$
Choose the first $h_j$ for which $f(s(p)+h_j)>0$ and put
$$
\widetilde p=s(p)+h_j.
$$
This is a measurable choice. Set $E=\Set{\widetilde p:p\in E_{\mathcal M}}$. Then $E\subseteq\Set{f>0}$. Therefore $E$ meets almost every $L_i$-coset exactly once, $i=1,2$.
\end{proof}

\section{The cases $k=0$ and $k=1$}\label{s:low-rank}

Let $\mathcal P$ be the convex set of measurable perfect fractional matchings on the bipartite graph $\mathcal{B}$. The weak compactness argument preceding \cite[Lemma 4.2]{bowen2021perfect} gives an extreme point $\sigma\in\mathcal P$. It is proved in \cite{bowen2021perfect} that
$$
\sigma(e)\in\Set{0,\frac12,1}
$$
for almost every edge, and that
$$
\mathcal L=\Set{e:\sigma(e)=\frac12}
$$
is, after deleting a null set, a vertex-disjoint union of bi-infinite paths. Every vertex outside $\mathcal L$ is already covered by a unique edge on which $\sigma=1$.

By Lemma \ref{lm:lift} it is enough to show that $\mathcal{B}$ has a measurable perfect matching, in order to show that some measurable $E \subseteq \Set{f > 0}$ tiles with both $L_1$ and $L_2$.

\subsection{The case $k=0$}

If $k=0$, the groups $\Lambda_i$ are finite. Hence every component of the support graphing is finite. A finite component cannot contain a bi-infinite path, so $\mathcal L$ is empty almost everywhere. Thus $\sigma$ is already a measurable perfect matching. 

\subsection{The case $k=1$}

Choose a surjective homomorphism
$$
q:\Lambda_2\longrightarrow\ZZ
$$
with finite kernel. The action $\Lambda_2\curvearrowright X_1$ is free by \eqref{eq:lambda-intersection}. Thus, for left vertices $x,x'$ in the same $\Lambda_2$-orbit, there is a unique $\gamma(x,x')\in\Lambda_2$ such that $x'=x+\gamma(x,x')$. Define
$$
c(x,x')=q(\gamma(x,x')).
$$
Then $c$ is Borel, integer-valued and satisfies
$$
c(x,z)=c(x,y)+c(y,z).
$$

If $x,x'$ are consecutive left vertices on a path of $\mathcal L$, their displacement $\gamma(x,x')$ belongs to the finite set $\Lambda_2\cap(K-K)$. Hence there is $C<\infty$ such that
$$
\Abs{c(x,x')}\leq C.
\eqn{eq:boundedstep}
$$
If $x_0$ is a left vertex on such a path and $h(x)=c(x_0,x)$, then
$$
\#\Set{x:h(x)=m}\leq\Abs{\ker q},
\qquad m\in\ZZ.
\eqn{eq:fibers}
$$

It follows from \eqref{eq:boundedstep} and \eqref{eq:fibers} that along either end of a path, $\Abs{h}$ tends to infinity and the sign of $h$ eventually stabilizes. The sign is independent of the choice of $x_0$. It is also Borel: a directed ray $x=x_0,x_1,\ldots$ is positive exactly when
$$
\text{for every }M\in\NN\text{ there is }J\text{ such that }c(x,x_j)\geq M\text{ for all }j\geq J.
$$

Let us show that the paths with two positive ends (similarly for two negative ends) form a null set. Indeed, on such a path $h$ tends to $+\infty$ at both ends and therefore has a nonempty finite set of global minima. Intrinsically, a left vertex $x$ is a minimum exactly when
$$
c(x,y)\geq0
\quad\text{for every left vertex }y\text{ on the same path}.
$$
The set of minima is Borel and contains between one and
$|\ker q|$ vertices of each such path. Hence one can select
Borel-measurably exactly one minimum from every path, obtaining
a Borel transversal $T$. Since the equivalence relation generated
by the graphing is probability-measure-preserving and each of
its classes is infinite, $T$ must be null.

Thus almost every path in $\mathcal L$ has one negative
and one positive end. Orienting each such path from
its negative to its positive end gives a measurable
perfect matching on $\mathcal L$, by choosing at each
left vertex the edge pointing toward the positive end.
The complementary edges give another measurable
perfect matching on $\mathcal L$. Extending both
matchings by the edges on which $\sigma=1$, we obtain
two measurable perfect matchings whose average is
$\sigma$. Since $\sigma$ is an extreme point of
$\mathcal P$, these matchings must agree almost
everywhere. Hence $\mathcal L$ is null and $\sigma$
itself is a measurable perfect matching.

\section{The case $k\geq2$}\label{s:higher-rank}

We use  \eqref{eq:convolution} to establish growth of the connected components on each side of the bipartite graph $\mathcal{B}$.

\begin{lemma}\label{lm:growth}
Let $D_1,D_2$ be bounded measurable fundamental domains of the lattices
$L_1,L_2$, respectively, and let $k=d-\operatorname{rank}(L_1\cap L_2) \ge 2$.
For almost every $x\in X_1$, the connected component of $x$ has superlinear growth.
\end{lemma}

\begin{proof}
Recall that $H=L_1\cap L_2$, $M=\RR^d/H$, and $\Lambda_i=L_i/H$.
We have $\Lambda_1\cap\Lambda_2=\{0\}$.

We first observe that the entire problem can be formulated on $M$. Indeed, for
$i=1,2$, define
$$
a_i(\pi(x))=\sum_{h\in H}\one_{D_i}(x+h).
$$
Since $H\subseteq L_i$ and $D_i$ is a fundamental domain for
$L_i$, the function $a_i$ takes only the values $0$ and $1$ almost everywhere.
Thus $a_i=\one_{A_i}$ almost everywhere for some bounded measurable set
$A_i\subseteq M$ which tiles $M$ by $\Lambda_i$.

Moreover, the quotient integration formula and the definition of $F$ give
\begin{equation}
F=\one_{A_1}*\one_{A_2}, \label{eq:quotient-convolution}
\end{equation}
where convolution is taken on $M$: 
expanding the convolution on $M$ and applying Weil's integration formula gives the periodization of $f$ over $H$.

\noindent \textbf{1. Definition of $T$.}

Let $Q$ be the closure of the set of density-one points of $A_2$.
Equivalently, $Q$ is the support of the measure $\one_{A_2}\,dm$.

Define $T$ to consist of those $\lambda\in\Lambda_2$ for which there exists
$z\in M$ such that every neighborhood of $z$ has positive measure intersection
with both $A_2$ and $A_2+\lambda$.

Equivalently,
$$
T=(Q-Q)\cap\Lambda_2.  
$$
Since $Q$ is compact and $\Lambda_2$ is discrete, $T$ is finite.
It is also symmetric.

\noindent \textbf{2. The set $T$ generates $\Lambda_2$.}

Since $A_2$ tiles $M$ by $\Lambda_2$, the translates of $Q$ by $\Lambda_2$
cover $M$. They form a locally finite family of closed sets.

Suppose that $\Gamma=\langle T\rangle$ is a proper subgroup of $\Lambda_2$.
Consider the sets
$$
U=\bigcup_{\gamma\in\Gamma}(Q+\gamma),
\qquad V=\bigcup_{\lambda\in\Lambda_2\setminus\Gamma} (Q+\lambda).
$$
Both are nonempty and closed, and together they cover $M$.

They are disjoint. Indeed, if $(Q+\gamma)\cap(Q+\lambda)\ne\varnothing$, then
$\lambda-\gamma\in T\subseteq\Gamma$, so $\gamma\in\Gamma$ implies
$\lambda\in\Gamma$.

Thus $M$ is the union of two disjoint nonempty closed sets, contradicting its
connectedness.

Consequently, $T$ generates $\Lambda_2$.

\noindent \textbf{3. The elements of T define a subgraph of the auxiliary graph $G$.}

Recall that $G$ is the auxiliary graph on $X_1=M/\Lambda_1$, in which two
distinct vertices are adjacent if they have a common neighbor in $X_2$.

We claim that, for every nonzero $\lambda\in T$, the vertices $x+\Lambda_1$ and
$x+\lambda+\Lambda_1$ are adjacent in $G$ for almost every $x\in M$.

Put $B_2=-A_2$. Since $T$ is symmetric and $\lambda\in T$, the essential
closures of $B_2$ and $B_2+\lambda$ intersect.  We can therefore find
density-one points $u_n\in B_2$ and $w_n\in B_2+\lambda$ such that
$$
\delta_n=w_n-u_n\longrightarrow0.
$$

Let $C_n\subseteq X_1$ consist of those classes $x+\Lambda_1$ for which $x+u_n$
and $x+w_n$ are density-one points of the same translate $A_1+a$, where
$a\in\Lambda_1$.

Since $A_1$ tiles $M$ by $\Lambda_1$, we have
$$
\mu_1(C_n) = m(A_1\cap(A_1-\delta_n)) \longrightarrow m(A_1)=1, 
$$
by translation continuity in $L^1(M)$.

Fix $x+\Lambda_1\in C_n$ and choose $a\in\Lambda_1$ as above.  The point
$x+u_n$ is a density-one point of both $A_1+a$ and $x+B_2$.  Consequently,
these two sets have positive-measure intersection.  Similarly, $x+w_n$ is a
density-one point of both $A_1+a$ and $x+B_2+\lambda$, so these sets also have
positive-measure intersection.

Using \eqref{eq:quotient-convolution}, we obtain
$$
F(x-a)>0,\qquad F(x+\lambda-a)>0. 
$$

Thus the two points
$$
p=x-a,\qquad p'=x+\lambda-a
$$
belong to $K=\{F>0\}$.

They determine edges of $\mathcal B$ whose left endpoints are $x+\Lambda_1$ and
$x+\lambda+\Lambda_1$, respectively.  Their right endpoints coincide, since
$p'-p=\lambda\in\Lambda_2$.

Hence these two left vertices are adjacent in $G$.

We have proved this adjacency for every $x+\Lambda_1\in C_n$. Since
$\mu_1(C_n)\to1$, the claim follows.

\noindent \textbf{4. Growth of the connected components.}

Since $T$ is finite and $\Lambda_2$ is countable, we can choose a
$\Lambda_2$-invariant conull subset $X_1^*\subseteq X_1$ such that, for every
$x\in X_1^*$ and every nonzero $\lambda\in T$, the vertices $x$ and $x+\lambda$
are adjacent in $G$.

Since $T$ generates $\Lambda_2$, all vertices in the $\Lambda_2$-orbit of $x\in
X_1^*$ belong to the same connected component of $G$.

Conversely, every edge of $G$ joins vertices whose difference belongs to
$\Lambda_2$. Thus the left vertices of almost every connected component of
$\mathcal B$ form exactly one $\Lambda_2$-orbit.

Now suppose $k\geq2$. Since $\Lambda_2\simeq\ZZ^k\times F_2$ and $T$ generates
$\Lambda_2$, we can choose two elements $\lambda_1,\lambda_2\in T$ which are
linearly independent in the free part of $\Lambda_2$.

For every $x\in X_1^*$, the vertices
$$
x+n_1\lambda_1+n_2\lambda_2, \qquad 0\leq n_1,n_2\leq N,
$$ are all
distinct, because $\Lambda_1\cap\Lambda_2=\{0\}$.

Each of these vertices can be reached from $x$ in at most $2N$ steps of $G$,
and hence in at most $4N$ steps of the bipartite graph $\mathcal B$.

Therefore
$$
|B_{\mathcal B}(x,4N)|\geq(N+1)^2.
$$
This proves quadratic growth for almost every left vertex. The same argument,
with the roles of the two lattices interchanged, gives the corresponding
estimate for almost every right vertex.

In particular, the connected components of $\mathcal B$ have superlinear growth
almost everywhere. 

\end{proof}

The bipartite graph $\mathcal{B}$ is hyperfinite by Lemma \ref{lm:graphing}. Since $k\geq2$, Lemma \ref{lm:growth} gives superlinear growth. Bowen, Kun and Sabok \cite[Proposition 2.1]{bowen2021perfect} show that a hyperfinite graphing with superlinear growth is almost everywhere one-ended. The fractional matching given by the periodization of $f$ is everywhere positive, so \cite[Theorem 1.1]{bowen2021perfect} gives a measurable perfect matching. Lemma \ref{lm:lift}, applied to $f$, produces a measurable common tile contained in
$$
\Set{f>0}\subseteq D_1+D_2.
$$
This proves Theorem \ref{th:main} in the remaining case $k\geq2$.

\section*{Acknowledgement.} We have used ChatGPT. It pointed us to using the result of Bowen, Kun and Sabok \cite{bowen2021perfect} and to how the case $k=1$ can be handled. All mathematical statements, proofs, references, and conclusions so obtained were reviewed
and verified by the authors, who take full responsibility for the content of the
paper.

\bibliographystyle{alpha}
\bibliography{mk-bibliography}

\end{document}